\documentclass[12pt,twoside]{amsart}
\usepackage{pstricks,pst-plot,}
\usepackage{amssymb}
\usepackage{tikz}
\usepackage{float}
\usepackage{verbatim}
\usepackage{amsmath}
\usepackage{bm}
\usepackage[margin=2.2cm, top=3.0cm, bottom=3.5cm]{geometry}
\usepackage[T1]{fontenc}
\usepackage{times}
\usepackage{amssymb,latexsym}
\usepackage{enumerate}
\usepackage[stable]{footmisc}
\usepackage{color}
\usepackage[colorlinks,linkcolor=black,citecolor=black,urlcolor=black]{hyperref}
\usepackage{amsthm}
\usepackage{mathrsfs}
\usepackage{geometry}
\usepackage{enumitem}
\usepackage{cases}
\usepackage{xcolor} 
\usepackage{appendix} 
\usepackage[utf8]{inputenc}

\usepackage{fancyhdr}
\usepackage{titlesec}
\usepackage{titletoc}
\titleformat{\section}
{\normalfont\large\bfseries}  
{\thesection}                 
{0.5em}                       
{}

\titleformat{\subsection}
{\normalfont\bfseries}        
{\thesubsection}              
{0.5em}                       
{}
\titleformat{\subsubsection}
{\normalfont\bfseries}        
{\thesubsubsection}              
{0.5em}                       
{}

\titlespacing*{\section}
{0pt}                         
{4ex plus 1ex minus 0.2ex}    
{2.5ex plus 0.2ex}            

\titlespacing*{\subsection}
{0pt}                         
{3.5ex plus 1ex minus 0.2ex}  
{2ex plus 0.2ex}              

\titlecontents{section}
[2em]
{\vspace{0.8ex}\normalfont}
{\contentslabel{2em}}
{\hspace*{-2em}}  
{\titlerule*[0.5pc]{.}\contentspage}

\titlecontents{subsection}
[3.5em]
{\vspace{0.5ex}\normalfont}
{\contentslabel{2.8em}}
{\hspace*{-2.8em}}
{\titlerule*[0.5pc]{.}\contentspage}

\titlecontents{subsection}
[3.5em]                        
{\vspace{0.5ex}\normalfont}    
{\contentslabel{2.8em}}        
{\hspace*{-2.8em}}             
{\titlerule*[0.5pc]{.}\contentspage}  

\makeatletter
\def\@afterheading{%
	\@nobreaktrue
	\everypar{%
		\if@nobreak
		\@nobreakfalse
		\clubpenalty\@M
		\else
		\clubpenalty\@clubpenalty
		\everypar{}%
		\fi}}
\makeatother

\newcommand{\bea}{\begin{eqnarray*}}
	\newcommand{\eea}{\end{eqnarray*}}

\makeatletter
\newcommand{\sumprime}{\if@display\sideset{}{'}\sum%
	\else\sum'\fi}
\makeatother

\numberwithin{equation}{section}

\newtheorem{theorem}{Theorem}
\newtheorem{lemma}{Lemma}
\newtheorem{proposition}{Proposition}

\newtheorem{example}{Example}
\newtheorem*{weightedrestatement}
  {Theorem~\ref{t031} (restated)}
  
\newcommand{\C}{\mathbb{C}}

\newcommand{\ddbar}{\bar\partial}

\title[Unique continuation property for $\bar\partial$ with $L^{2n}$ potentials]{Unique continuation property for $\bar\partial$ operator with $L^{2n}$ potentials}

\author{Guokuan  Shao}
\address{School of Mathematics (Zhuhai), Sun Yat-Sen University, Zhuhai, Guangdong 519082, China}
\email{shaogk@mail.sysu.edu.cn}
\thanks{Guokuan Shao was supported by National Key R\& D Program of China No. 2024YFA1015200, National Natural Science Foundation of China No. 12471082. Ziming Shi was supported by SUSTech research start-up grant Y01286157 and Y01286257. Jujie Wu was supported by National Key R\&D Program of China, No. 2024YFA1015200, National Natural Science Foundation of China,  No. 12671104.} 

\author{Ziming Shi}
\address{School of Mathematical Sciences, Southern University of Science and Technology, Shenzhen, Guangdong 518055, China} 
\email{shizm@sustech.edu.cn}

\author{Jujie Wu}
\address{School of Mathematics (Zhuhai), Sun Yat-Sen University, Zhuhai, Guangdong 519082, China}
\email{wujj86@mail.sysu.edu.cn}

\begin{document}

\begin{abstract}
	We establish the sharp unique continuation property for the Cauchy--Riemann operator $\bar{\partial}$ in higher dimensions. Let $0\in\Omega \subset \C^n$ be a domain. Suppose $u \in W_{\mathrm{loc}}^{1,2}(\Omega; \C^N)$ satisfies the differential inequality
	\[
	|\bar{\partial} u| \leq V |u| \quad \text{a.e.~on } \Omega
	\]
	with $V \in L_{\mathrm{loc}}^{2n}(\Omega)$. If $u$ is $L^2$ flat  at the origin, then $u \equiv 0$ on $\Omega$. 
	
	\bigskip

	
\end{abstract}

\maketitle
\tableofcontents  

\section{Introduction}

The unique continuation property (UCP) for partial differential equations concerns a fundamental question in the theory of PDEs: if a solution of a differential equation or inequality vanishes to infinite order at a point, or vanishes on an open subset of its domain, does it vanish identically throughout the connected component of the domain? For elliptic equations of second order with sufficiently regular coefficients, classical results of Carleman \cite{C}, Aronszajn \cite{A},  H\"ormander \cite{H} and Jerison-Kenig \cite{JK} established weak and strong UCP under appropriate integrability assumptions. 
In this paper, we study strong unique continuation
for the first-order Cauchy--Riemann operator
$\bar\partial$ under critical integrability
assumptions on the potential.

Let $\Omega \subset \C^n$ be a domain and consider a vector-valued function $u = (u_1, \ldots, u_N): \Omega \to \C^N$. We say that $u$ satisfies the differential inequality
\begin{equation}\label{e83106}
|\bar{\partial} u| \leq V |u| \quad \text{a.e.~on } \Omega, 
\end{equation}
where $V  \geq 0 $ is a measurable function (the ``potential''). The above inequality encompasses a wide range of equations, including the homogeneous $\bar{\partial}$-equation (when $V \equiv 0$), linear systems and more general Schr\"odinger-type systems.

A function $u \in L^2_{\mathrm{loc}}(\Omega, \C^N)$ is said to be $L^2$-flat  at a point $z_0 \in \Omega$ if, for every $M \geq 0$,
\begin{equation}\label{e83107}
\lim_{r \to 0} r^{-M} \int_{|z - z_0| < r} |u(z)|^2 \, dV = 0. 
\end{equation}
The \emph{strong} UCP holds if every solution of \eqref{e83106} that is $L^2$ flat at some point vanishes identically. The \emph{weak} UCP holds if every solution that vanishes on a nonempty open subset vanishes identically.
 
Throughout this paper, statements that a Sobolev function
vanishes are understood in the almost-everywhere sense.

The study of UCP for $\bar{\partial}$-type equations has a long history, motivated by problems in several complex variables, CR-geometry. Bell-Lempert \cite{BellLempert} proved and applied the strong UCP for $|\bar{\partial} u| \leq C |u|$ with a constant potential $C$ while studying boundary regularity and uniqueness for CR-mappings of hypersurfaces. Their result raised the natural question: what is the minimal integrability assumption on the potential $V$ that still guarantees UCP?

Denote by ${\rm W}^{k,2}_{\mathrm{loc}}(\Omega)$ the standard (local) Sobolev space
of functions whose weak derivatives up to order $k$ are in $L^2_{\rm{loc}}(\Omega)$.
For the one-dimensional case $n = 1$, the equation $\bar{\partial} u = V u$ can be studied via the similarity principle and the theory of generalized analytic functions when $V\in L^p$, $p>2$.  The critical integrability threshold is $L^2_{\mathrm{loc}}$.
Indeed, Pan and Zhang \cite[Theorem~1.1]{PanZhangL2} proved that,
if $\Omega\subset\mathbb C$ is a domain and
$u\in {\rm W}^{1,2}_{\mathrm{loc}}(\Omega;\mathbb C^N)$ satisfies
$|\bar\partial u|\le V|u|$ a.e. on $\Omega$, with
$V\in L^2_{\mathrm{loc}}(\Omega)$, then vanishing to
infinite order at a point in the $L^2$ sense implies
that $u=0$ a.e. on $\Omega$. 
This result is sharp: Mandache \cite{Mandache} constructed counterexamples showing that the weak UCP fails for potentials in $L^p_{\mathrm{loc}}$ with $0 < p < 2$. 

In higher dimensions $n \geq 2$, the situation is substantially more subtle. The natural integrability threshold for strong UCP should be $L^{2n}_{\mathrm{loc}}$, as  indicated  by Pan-Zhang \cite{PanZhangL2} and Shi \cite{Shi}. However, for weak UCP, a surprising phenomenon occurs: the threshold remains $L^2_{\mathrm{loc}}$, independent of the dimension $n$. This was established by Pan-Zhang \cite{PanZhangL2}, where they proved that weak UCP holds for (1.1) with $V \in L^2_{\mathrm{loc}}(\Omega)$ in all dimensions $n \geq 2$. Their proof proceeds by slicing along coordinate axes and applying the one-dimensional strong UCP on almost every slice, exploiting Fubini's theorem to preserve the $L^2$ integrability of the potential.

The present paper addresses the sharp strong UCP for $\bar{\partial}$ in higher dimensions. 
The main result of this paper establishes strong
unique continuation for vector-valued
$ {\rm W}^{1,2}_{\rm loc}$ solutions of
$|\bar\partial u|\le V|u|$ with
$V\in L^{2n}_{\mathrm{loc}}$, in every complex
dimension $n\ge1$ and for every target dimension
$N\ge1$.
\begin{theorem}\label{t01}
Let $n,N\ge1$ be integers, and let
$\Omega\subset\mathbb C^n$ be a domain containing the origin.
Suppose that
$u\in W^{1,2}_{\mathrm{loc}}(\Omega;\mathbb C^N)$
satisfies
\[
|\bar\partial u|\le V|u|
\quad\text{a.e. on }\Omega,
\]
where $V\ge0$ belongs to $L^{2n}_{\mathrm{loc}}(\Omega)$.
 If $u$ is $L^2$-flat at the origin,  then $u$  vanishes identically on $\Omega$.
\end{theorem}
The exponent $2n$ is optimal within the Lebesgue scale. For every $0< p < 2n$, there exist nontrivial smooth functions $u$ that are flat and satisfy $|\bar{\partial} u| \leq V |u|$ with $V \in L^p_{\mathrm{loc}}$ but $V \notin L^{2n}_{\mathrm{loc}}$ (cf. Example \ref{ex1}). 

Once the reduction to one dimension is achieved, we apply a new one-dimensional critical weighted UCP as follows. 

\begin{theorem}\label{t031}
Let $n,N\ge1$ be integers and set $s=(n-1)/n$.
Let $\Omega\subset\mathbb C$ be a domain containing
the origin. Suppose that
$v\in {\rm W}^{1,2}_{\mathrm{loc}}(\Omega;\mathbb C^N)$
satisfies
\[
|\bar\partial v(z)|
\le |z|^{-s}W(z)|v(z)|
\quad\text{a.e. on }\Omega\setminus\{0\},
\]
where $W\ge0$ belongs to
$L^{2n}_{\mathrm{loc}}(\Omega)$.
If $v$ is   $L^2$-flat at the origin,  then  $v$  vanishes identically on $\Omega$.
\end{theorem}


For $n\ge3$ and $N\ge2$, Theorem~\ref{t031}
gives an affirmative answer to the weighted planar
unique continuation question posed by Pan and Zhang
\cite[Remark~7.2, Question~2]{PZ2025}
for smooth solutions, as well as to its
${\rm W}^{1,2}_{\mathrm{loc}}$ formulation in
\cite[Remark~4.2, Question~2]{PanZhang2}.
It will be proved via a carefully designed Carleman inequality with optimal exponents. Specifically, we establish that for all $f \in C_c^\infty(\C\setminus\{0\}; \C^N)$,
\[
\||z|^{-\tau_m} f\|_{L^{q_n}(\C)} \leq C_n \||z|^{-\tau_m + s} \bar{\partial} f\|_{L^{p_n}(\C)}
\]
for suitable parameters and exponents.
The proof of this one-dimensional inequality relies on the Stein-Weiss theorem \cite{SW} on weighted Riesz potentials, applied with the homogeneity parameters chosen precisely to match the exponents. 
One key observation is that
\[
    p_n=\frac{4n}{2n+1},\qquad q_n=\frac{4n}{2n-1},\qquad
    \frac1{p_n}-\frac1{q_n}=\frac1{2n}.
\] 
The identity $1/p_n=1/(2n)+1/q_n$ yields
\[
\|WF\|_{L^{p_n}}
\le \|W\|_{L^{2n}}\|F\|_{L^{q_n}}.
\]
The identity $1/p_n=1/(2n)+1/q_n$ allows us to estimate
the potential term by H\"older's inequality. By the absolute
continuity of the $L^{2n}$ integral, this term can be absorbed
on a sufficiently small disk. The exponent $2n$ is also
critical with respect to the scaling of the weighted
differential inequality, as explained in Section~2.

The compensated slicing technique developed here has several important consequences.  First, it establishes strong unique continuation at the critical Lebesgue exponent for
${\rm W}^{1,2}_{\rm loc}$ vector-valued solutions
in arbitrary complex dimension. 
  Second, our finite-point strong unique continuation
results complement the work of Pan and Zhang
\cite{PZ} on unique continuation at infinity. 
Third, the technique is likely to have broader applications to other first-order elliptic systems and to problems involving complex geometry where polar weights and holomorphic factors interact.

The paper is organized as follows. In Section 2, we prove the one-dimensional critical weighted UCP, establishing the key Carleman inequality. In Section 3, we develop the compensated slicing technique and prove the main theorem. Appendix~A establishes the failure, for $n\ge2$,
of the uniform radial $L^{p_n}$--$L^{q_n}$ Carleman
estimate considered in Section~3. This rules out
a direct proof based on that estimate and motivates
our reduction to a one-dimensional weighted problem
through compensated slicing.


\section{Critical weighted strong unique continuation in one complex variable}
\label{sec:weighted-ucp}

We establish strong unique continuation for a planar
$\bar\partial$ inequality with a singular weight at the origin.
Throughout this section, $n,N\geq1$ are integers and
\begin{equation}\label{eq:ucp-parameters}
    s=\frac{n-1}{n},\qquad a=\frac{s}{2},\qquad
    p_n=\frac{4n}{2n+1},\qquad q_n=\frac{4n}{2n-1}.
\end{equation}
Here $n$ parametrizes the weight and the integrability exponent;
the source domain is always contained in $\mathbb C$.
We write $D_r=\{z\in\mathbb C:|z|<r\}$, denote planar Lebesgue
measure by $dA$, and use the Euclidean norm on $\mathbb C^N$.
The operator $\bar\partial$ denotes
$\partial/\partial\bar z=\tfrac12(\partial_x+i\partial_y)$,
acting componentwise. 

\begin{weightedrestatement}
Let $\Omega\subset\mathbb C$ be a domain containing the origin,
and let $W\geq0$ belong to $L^{2n}_{\mathrm{loc}}(\Omega)$.
Suppose that $v\in {\rm W}^{1,2}_{\mathrm{loc}}(\Omega;\mathbb C^N)$
satisfies
\begin{equation}\label{eq:weighted-dbar-inequality}
    |\bar\partial v(z)|\leq |z|^{-s}W(z)|v(z)|
    \qquad\text{for almost every }z\in\Omega\setminus\{0\}.
\end{equation}
If $v$ vanishes to infinite order at the origin in the $L^2$ sense,
that is,
\begin{equation}\label{eq:l2-flatness}
    \lim_{r\rightarrow 0}r^{-M}\int_{D_r}|v|^2\,dA=0
    \qquad\text{for every }M>0,
\end{equation}
then $v\equiv0$ on $\Omega$.
\end{weightedrestatement}

The exponent $2n$ is critical with respect to dilation:
if $v_t(z)=v(tz)$ for $t>0$, then the rescaled inequality has
coefficient
\[
    W_t(z)=t^{1-s}W(tz),\qquad
    \|W_t\|_{L^{2n}(U)}=\|W\|_{L^{2n}(tU)},
\]
for every measurable $U\subset t^{-1}\Omega$, because
$1-s=1/n$. The proof combines a weighted Carleman estimate
with the local smallness of the $L^{2n}$ norm of $W$.

\subsection{The weighted Carleman estimate}

We use the following form of the Stein--Weiss inequality
\cite{SW}. For $d\geq1$ and $0<\alpha<d$, define
\[
    I_\alpha h(x)
    =\int_{\mathbb R^d}\frac{h(y)}{|x-y|^{d-\alpha}}\,dy.
\]

\begin{theorem}[Stein--Weiss, \cite{SW} ]\label{t02}
Let $1<p\leq q<\infty$ and $p'=p/(p-1)$.
Suppose that $\beta,\gamma\in\mathbb R$ satisfy
\[
    \beta<\frac{d}{q},\qquad
    \gamma<\frac{d}{p'},\qquad
    \beta+\gamma\geq0,\qquad
    \frac1p-\frac1q=\frac{\alpha-\beta-\gamma}{d}.
\]
Then
\[
    \bigl\||x|^{-\beta}I_\alpha h\bigr\|_{L^q(\mathbb R^d)}
    \leq C\bigl\||x|^\gamma h\bigr\|_{L^p(\mathbb R^d)}
\]
for every measurable $h$ with $|x|^\gamma h\in L^p(\mathbb R^d)$,
where $C$ depends only on $d,\alpha,p,q,\beta,\gamma$.
\end{theorem}

The homogeneity condition can be checked by setting $h_t(x)=h(tx)$.
Indeed, $I_\alpha h_t(x)=t^{-\alpha}(I_\alpha h)(tx)$, and the
two norms scale by $t^{\beta-\alpha-d/q}$ and
$t^{-\gamma-d/p}$, respectively. Applying the inequality to
all dilates of a fixed nonzero, nonnegative test function forces
these powers to agree.

Let $s, a, p_n, q_n$ be as in \eqref{eq:ucp-parameters}.
Then we have
\begin{equation}\label{eq:exponent-identities}
    1<p_n<2<q_n<\infty,\qquad
    p_n'=q_n,\qquad
    \frac1{p_n}=\frac1{2n}+\frac1{q_n}.
\end{equation}
Moreover,
\[
    \frac{2}{q_n}-a=\frac{2}{p_n'}-a=\frac12>0,
    \qquad 2a=s\geq0,
    \qquad
    \frac1{p_n}-\frac1{q_n}=\frac{1-2a}{2}.
\]
Thus, Theorem~\ref{t02}, with $d=2$, $\alpha=1$, and
$\beta=\gamma=a$, yields
\begin{equation}\label{e82908}
    \bigl\||z|^{-a}I_1h\bigr\|_{L^{q_n}(\mathbb C)}
    \leq A_n\bigl\||z|^a h\bigr\|_{L^{p_n}(\mathbb C)},
    \qquad
    I_1h(z)=\int_{\mathbb C}\frac{h(\zeta)}{|z-\zeta|}\,dA(\zeta),
\end{equation}
where $A_n$ depends only on $n$.

\begin{proposition} \label{p01}
Set $\tau_m=m+a$ for $m=0,1,2,\ldots$.
There exists $C_n>0$, depending only on $n$, such that
\begin{equation}\label{e90401}
    \bigl\||z|^{-\tau_m}f\bigr\|_{L^{q_n}(\mathbb C)}
    \leq C_n
    \bigl\||z|^{-\tau_m+s}\bar\partial f
    \bigr\|_{L^{p_n}(\mathbb C)}
\end{equation}
for every $f\in C_c^\infty(\mathbb C\setminus\{0\};\mathbb C^N)$.
The same estimate holds for every
$f\in W^{1,p_n}(\mathbb C;\mathbb C^N)$ whose support is a compact
subset of $\mathbb C\setminus\{0\}$.
\end{proposition}

\begin{proof}
For $g\in C_c^\infty(\mathbb C;\mathbb C^N)$, the componentwise
Cauchy--Pompeiu formula and the triangle inequality in $\mathbb C^N$
give
\[
    g(z)=\frac1\pi\int_{\mathbb C}
      \frac{\bar\partial g(\zeta)}{z-\zeta}\,dA(\zeta),
    \qquad
    |g(z)|\leq\frac1\pi I_1(|\bar\partial g|)(z).
\]
Applying \eqref{e82908} to the scalar function
$|\bar\partial g|$, we obtain
\begin{equation}\label{e82909}
    \bigl\||z|^{-a}g\bigr\|_{L^{q_n}}
    \leq C_n\bigl\||z|^a\bar\partial g\bigr\|_{L^{p_n}},
    \qquad C_n=A_n/\pi.
\end{equation}
This also proves that the constant is independent of $N$.

Given $f\in C_c^\infty(\mathbb C\setminus\{0\};\mathbb C^N)$,
set $g(z)=z^{-m}f(z)$. Since $f$ vanishes near the origin, $g$
extends by zero to a smooth compactly supported function on
$\mathbb C$. Holomorphicity of $z^{-m}$ away from the origin gives
\[
    \bar\partial g=z^{-m}\bar\partial f.
\]
Substituting into \eqref{e82909} and using
$-m-a=-\tau_m$ and $-m+a=-\tau_m+s$ proves \eqref{e90401},
with the same constant for every $m$.

To obtain the Sobolev extension, fix $m$ and let
$f\in W^{1,p_n}(\mathbb C;\mathbb C^N)$ have compact support away
from the origin. Its standard mollifications $f_j$ converge to $f$
in $W^{1,p_n}(\mathbb C)$ and, after passage to a subsequence,
almost everywhere. For all sufficiently large $j$, their supports
lie in a fixed compact annulus. Both weights in \eqref{e90401}
are bounded on this annulus. Therefore
\[
    |z|^{-\tau_m+s}\bar\partial f_j
    \longrightarrow |z|^{-\tau_m+s}\bar\partial f
    \quad\text{in }L^{p_n}(\mathbb C).
\]
Applying \eqref{e90401} to $f_j$ and then using Fatou's lemma yields
\[
\begin{aligned}
    \bigl\||z|^{-\tau_m}f\bigr\|_{L^{q_n}}
    &\leq\liminf_{j\to\infty}
      \bigl\||z|^{-\tau_m}f_j\bigr\|_{L^{q_n}}\\
    &\leq C_n
      \bigl\||z|^{-\tau_m+s}\bar\partial f\bigr\|_{L^{p_n}}.
\end{aligned}
\]
In particular, the weighted $L^{q_n}$ norm on the left is finite.
\end{proof}

\subsection{Proof of  Theorem ~\ref{t031}}

\begin{proof}[Proof of Theorem~\ref{t031}]
Write now $p=p_n$ and $q=q_n$.
By absolute continuity of the integral, we may choose $R>0$ so that
\begin{equation}\label{eq:local-smallness}
    D_{2R}\Subset\Omega,\qquad
    C_n\|W\|_{L^{2n}(D_{2R})}\leq\frac12,
\end{equation}
where $C_n$ is the constant in Proposition~\ref{p01}.
This choice of $R$ is independent of $m$.

\emph{Step 1: Localization and absorption.}
Choose $\eta\in C_c^\infty(D_{2R})$ with
$0\leq\eta\leq1$ and $\eta=1$ on $D_R$.
For $0<\varepsilon<R/4$, choose a smooth radial cutoff
$\chi_\varepsilon$ such that
\[
    0\leq\chi_\varepsilon\leq1,\qquad
    \chi_\varepsilon=0\ \text{on }D_\varepsilon,\qquad
    \chi_\varepsilon=1\ \text{outside }D_{2\varepsilon},\qquad
    |\bar\partial\chi_\varepsilon|\leq C/\varepsilon.
\]
The function $f_\varepsilon=\eta\chi_\varepsilon v$, extended
by zero outside $\Omega$, belongs to $W^{1,2}(\mathbb C;\mathbb C^N)$
and has compact support away from the origin. Since $p<2$,
it also belongs to $W^{1,p}(\mathbb C;\mathbb C^N)$.
By Proposition~\ref{p01}  we have
 $$
\bigl\||z|^{-\tau}f_\varepsilon\bigr\|_{L^q(\C)} \leq C_n
    \bigl\||z|^{-\tau+s}\bar\partial f_\varepsilon
    \bigr\|_{L^{p_n}(\mathbb C)}.
$$
To estimate the right-hand side, note that $\overline{\partial} f _\varepsilon = \eta
               (\bar\partial\chi_\varepsilon)v+ (\bar\partial\eta)\chi_\varepsilon v+ \eta \chi_\varepsilon \overline{\partial} v.$
We may fix $m\geq0$, write $\tau=\tau_m$, and set
\[
\begin{aligned}
    X_{\varepsilon,m}
      &=\bigl\||z|^{-\tau}f_\varepsilon\bigr\|_{L^q},\\
    E_{\varepsilon,m}
      &=\bigl\||z|^{-\tau+s}\eta
               (\bar\partial\chi_\varepsilon)v\bigr\|_{L^p},\\
    F_m
      &=\bigl\||z|^{-\tau+s}(\bar\partial\eta)v\bigr\|_{L^p}.
\end{aligned}
\]
Then the triangle inequality implies 
\[
    X_{\varepsilon,m}
    \leq C_n\bigl\||z|^{-\tau+s}\eta\chi_\varepsilon
                     \bar\partial v\bigr\|_{L^p}
       +C_n(E_{\varepsilon,m}+F_m).
\]
The differential inequality and H\"older's inequality, with
$1/p=1/(2n)+1/q$, give
\[
\begin{aligned}
    \bigl\||z|^{-\tau+s}\eta\chi_\varepsilon
                     \bar\partial v\bigr\|_{L^p}
    &\leq\bigl\|W|z|^{-\tau}f_\varepsilon\bigr\|_{L^p}\\
    &\leq\|W\|_{L^{2n}(D_{2R})}X_{\varepsilon,m}.
\end{aligned}
\]
Since $X_{\varepsilon,m}<\infty$, \eqref{eq:local-smallness}
allows this term to be absorbed, yielding
\begin{equation}\label{eq:absorbed-carleman}
    X_{\varepsilon,m}\leq2C_n(E_{\varepsilon,m}+F_m).
\end{equation}

\emph{Step 2: Removal of the inner cutoff.}
The derivative of $\chi_\varepsilon$ is supported in
$D_{2\varepsilon}\setminus D_\varepsilon$.
Since $-\tau+s=-m+a$ and $m\geq0$, on this annulus we have
$|z|^{-\tau+s}\leq2^a\varepsilon^{-m+a}$.
Consequently,
\begin{equation}\label{eq:inner-cutoff-error}
\begin{aligned}
    E_{\varepsilon,m}
    &\leq C\varepsilon^{-m+a-1}
                \|v\|_{L^p(D_{2\varepsilon})}\\
    &\leq C\varepsilon^{-m+a-1+2(1/p-1/2)}
                \|v\|_{L^2(D_{2\varepsilon})}\\
    &=C\varepsilon^{-m-1/2}
                \|v\|_{L^2(D_{2\varepsilon})}.
\end{aligned}
\end{equation}
Here we used $2(1/p-1/2)=1/(2n)$ and $a+1/(2n)=1/2$.
For each fixed $m$, the flatness assumption implies
\[
    \varepsilon^{-2m-1}
    \int_{D_{2\varepsilon}}|v|^2\,dA\longrightarrow0,
\]
and hence $E_{\varepsilon,m}\to0$ as $\varepsilon\rightarrow 0$.
As $f_\varepsilon\to\eta v$ almost everywhere,
Fatou's lemma and \eqref{eq:absorbed-carleman} yield
\begin{equation}\label{eq:after-inner-cutoff}
    \bigl\||z|^{-\tau_m}\eta v\bigr\|_{L^q}
    \leq2C_nF_m.
\end{equation}

\emph{Step 3: Vanishing near the origin.}
The support of $\bar\partial\eta$ is contained in
$\{R\leq|z|<2R\}$. Since $\tau_m\geq0$,
\[
    F_m\leq R^{-\tau_m}M_R,\qquad
    M_R:=\bigl\||z|^s(\bar\partial\eta)v\bigr\|_{L^p}<\infty.
\]
The quantity $M_R$ is independent of $m$.
For any $0<\rho<R$, using $\eta=1$ on $D_R$ in
\eqref{eq:after-inner-cutoff}, we obtain
\[
    \rho^{-\tau_m}\|v\|_{L^q(D_\rho)}
    \leq\bigl\||z|^{-\tau_m}\eta v\bigr\|_{L^q}
    \leq2C_nR^{-\tau_m}M_R.
\]
It follows that
\begin{equation}\label{eq:geometric-decay}
    \|v\|_{L^q(D_\rho)}
    \leq2C_nM_R\left(\frac{\rho}{R}\right)^{\tau_m}.
\end{equation}
Now let $m\to\infty$. Since $\tau_m\to\infty$ and $\rho<R$,
we conclude that $v=0$ on $D_\rho$, and hence on $D_R$.
Thus the limits have been taken first in $\varepsilon$ for each
fixed $m$, and only afterwards in $m$.

\emph{Step 4: Propagation to the whole domain.}
On the connected planar domain $\Omega^*=\Omega\setminus\{0\}$,
the coefficient $V(z)=|z|^{-s}W(z)$ belongs to
$L^2_{\mathrm{loc}}(\Omega^*)$: the factor $|z|^{-s}$ is locally
bounded there and $2n\geq2$. The weak unique continuation theorem
for $\bar\partial$ with square-integrable potentials
\cite[Theorem~1.3]{PanZhangL2} therefore applies to $v$ on
$\Omega^*$. Since $v$ vanishes on the nonempty open subset
$D_R\setminus\{0\}$, it vanishes throughout $\Omega^*$.
Together with the local conclusion at the origin, this proves
$v\equiv0$ on $\Omega$.
\end{proof}

\section{Proof of Theorem \ref{t01}}
In this section we give the proof of the main theorem. The main difficulty is that the natural Carleman inequality similar to the inequality for Laplacian in \cite{JK}  does not hold for $\bar\partial$ operator. 

Let $p_n=\frac{4n}{2n+1}$, $q_n=\frac{4n}{2n-1}$, where $n\geq 2$. Suppose that there exists $C>0$ independent of $\tau$ such that
\begin{equation}\label{e01}
	\||z|^{-\tau} f\|_{L^{q_n}(\mathbb{C}^n)} \le C \||z|^{-\tau} \bar\partial f\|_{L^{p_n}(\mathbb{C}^n)} 
\end{equation}
for all $f \in C_c^\infty(\mathbb{C}^n \setminus \{0\},\C^N)$. 
Such a uniform estimate would imply strong unique
continuation by a standard localization and
absorption argument.
For $n\ge2$, however, this estimate fails to hold
with a constant independent of $\tau$, even when
$\tau$ is restricted to an arbitrary sequence
tending to infinity; see Appendix~A.
This rules out a direct proof based on the above
uniform estimate and motivates a reduction to
the one-dimensional weighted setting.


We apply a \emph{compensated slicing technique} inspired by Pan-Zhang's work (cf. \cite{PZ2025}) that overcomes the failure of the Carleman inequality. 

Fix $0<r<R$ such that $B_R\Subset\Omega$.
For $\zeta\in S^{2n-1}$, define
\[
g_\zeta(w):=w^{n-1}u(w\zeta),
\qquad w\in D_r.
\]
The holomorphic factor $w^{n-1}$ serves three essential purposes.
First it exactly absorbs the polar Jacobian $|w|^{2n-2}$ in the complex polar identity \cite{PZ2025}
	\begin{eqnarray} \label{e82920}
	\int_{S^{2n-1}} \int_{D_r} |w|^{2n-2} F(w\zeta) \, dA(w) \, d\sigma(\zeta) = 2\pi \int_{B_r} F(z) \, dV(z).
	\end{eqnarray}
Here $d\sigma$ denotes the unnormalized Euclidean
surface measure on $S^{2n-1}$.
This identity provides $L^2$ control of the compensated
slices. Full $ {\rm W}^{1,2}_{\rm loc}$  control additionally requires a local
Hardy inequality to estimate the term
$(n-1)w^{n-2}u(w\zeta)$ arising from differentiation
of the holomorphic factor with respect to $w$.
Combining these estimates yields
$g_\zeta\in W^{1,2}(D_r;\mathbb C^N)$
for almost every $\zeta\in S^{2n-1}$,
as established in the compensated slicing lemma below.    
Second it preserves the $\bar{\partial}$-inequality in the form
	\[
	|\bar{\partial}_w g_\zeta(w)| \leq |w|^{-s} W_\zeta(w) |g_\zeta(w)|, 
	\]
	where $s = \frac{n-1}{n}$ and $W_\zeta(w) := |w|^s V(w\zeta) \in L^{2n}(D_r)$ for almost every $\zeta$.
Third it transfers the $L^2$-flatness of $u$ at the origin to $L^2$-flatness of $g_\zeta$ at $w = 0$ for almost every $\zeta$, through a dyadic argument using the Borel-Cantelli lemma.




For a smooth map $h$, define
\begin{equation}\label{e82922}
T_\zeta h(w) := w^{n-1} h(w \zeta), \quad   w\in D_r, \quad \zeta \in S^{2n-1}.
\end{equation}
The following compensated slicing lemma is crucial in the proof of the main theorem.

\begin{lemma}\label{l01}
With the notation above, for every smooth map $h$, we have
\[
\int_{S^{2n-1}}
\|T_\zeta h\|_{{\rm W}^{1,2}(D_r)}^2\,d\sigma(\zeta)
\le C_{n,r,R}\|h\|_{{\rm W}^{1,2}(B_R)}^2.
\]
Moreover, if
$u\in {\rm W}^{1,2}_{\rm loc}(\Omega;\mathbb C^N)$,
then
\[
T_\zeta u(w):=w^{n-1}u(w\zeta)
\]
belongs to ${\rm W}^{1,2}(D_r;\mathbb C^N)$ for almost
every $\zeta\in S^{2n-1}$, and the same estimate holds
with $h$ replaced by $u$.
\end{lemma}

\begin{proof}
We first consider the smooth case. We have
\begin{equation}\label{e83001}
\bar\partial_w T_\zeta h(w)
=
w^{n-1}\sum_{j=1}^n
\bar\zeta_j\bar\partial_{z_j}h(w\zeta),
\end{equation}
because $w^{n-1}$ is holomorphic. Similarly,
\begin{equation}\label{e83002}
\partial_w T_\zeta h(w)
=
(n-1)w^{n-2}h(w\zeta)
+
w^{n-1}\sum_{j=1}^n
\zeta_j\partial_{z_j}h(w\zeta).
\end{equation}
When $n=1$, the first term on the right is understood
to be zero.

By \eqref{e82920}, we have
\begin{equation}\label{e83002d}
\begin{split}
&\int_{S^{2n-1}}\int_{D_r}
|T_\zeta h(w)|^2\,dA(w)\,d\sigma(\zeta)\\
&\qquad=
\int_{S^{2n-1}}\int_{D_r}
|w|^{2n-2}|h(w\zeta)|^2\,dA(w)\,d\sigma(\zeta)\\
&\qquad=
2\pi\int_{B_r}|h(z)|^2\,dV(z).
\end{split}
\end{equation}

Since $|\zeta|=1$,
\[
\left|
\sum_{j=1}^n\bar\zeta_j\bar\partial_{z_j}h
\right|
\le |\bar\partial h|.
\]
Consequently,
\begin{equation}\label{e83003}
\int_{S^{2n-1}}\int_{D_r}
|\bar\partial_w T_\zeta h|^2\,dA\,d\sigma
\le
2\pi\int_{B_r}|\bar\partial h|^2\,dV.
\end{equation}
From \eqref{e83002}, we obtain
\begin{equation}\label{e83004}
|\partial_w T_\zeta h|^2
\le
2(n-1)^2|w|^{2n-4}|h(w\zeta)|^2
+
2|w|^{2n-2}|\partial h(w\zeta)|^2.
\end{equation}
The integral of the second term in \eqref{e83004}
is controlled by $\|\partial h\|_{L^2(B_r)}^2$
through \eqref{e82920}.
The first term vanishes when $n=1$.
For $n\ge2$, the local Hardy inequality
(cf.\ Appendix B) gives
\begin{equation}\label{e83005}
\begin{split}
&\int_{S^{2n-1}}\int_{D_r}
|w|^{2n-4}|h(w\zeta)|^2\,dA\,d\sigma\\
&\qquad=
2\pi\int_{B_r}\frac{|h(z)|^2}{|z|^2}\,dV(z)\\
&\qquad\le
C_{n,r,R}\|h\|_{{\rm W}^{1,2}(B_R)}^2,
\end{split}
\end{equation}
where $C_{n,r,R}$ depends only on $n,r,R$.
Combining \eqref{e83002d}, \eqref{e83003},
\eqref{e83004}, and, when $n\ge2$, \eqref{e83005},
we obtain
\begin{equation}\label{e83006}
\int_{S^{2n-1}}
\|T_\zeta h\|_{{\rm W}^{1,2}(D_r)}^2\,d\sigma(\zeta)
\le
C_{n,r,R}\|h\|_{{\rm W}^{1,2}(B_R)}^2.
\end{equation}

We now consider
$u\in {\rm W}^{1,2}_{\rm loc}(\Omega;\mathbb C^N)$.
Since $B_R\Subset\Omega$, we have
$u\in {\rm W}^{1,2}(B_R;\mathbb C^N)$.
Let
\[
\mathbb X:={\rm W}^{1,2}(D_r;\mathbb C^N).
\]
Initially, for smooth $h$, define
\[
(Th)(\zeta)(w):=T_\zeta h(w)=w^{n-1}h(w\zeta).
\]
We can rewrite \eqref{e83006} as
\begin{equation}\label{e83007}
\|Th\|_{L^2(S^{2n-1};\mathbb X)}
\le
C_{n,r,R}^{1/2}\|h\|_{{\rm W}^{1,2}(B_R)}.
\end{equation}
Thus $T$ is bounded with respect to the
${\rm W}^{1,2}(B_R)$ norm on smooth functions
and the norm of the Bochner space
$L^2(S^{2n-1};\mathbb X)$.
By density of smooth functions in
${\rm W}^{1,2}(B_R;\mathbb C^N)$,
$T$ extends uniquely to a bounded linear operator
\[
T:
{\rm W}^{1,2}(B_R;\mathbb C^N)
\longrightarrow
L^2\bigl(S^{2n-1};
{\rm W}^{1,2}(D_r;\mathbb C^N)\bigr).
\]
In particular, $(Tu)(\zeta)$ belongs to
${\rm W}^{1,2}(D_r;\mathbb C^N)$ for almost every
$\zeta$.
It remains to identify $(Tu)(\zeta)$ with the
compensated slice $w^{n-1}u(w\zeta)$.

Choose smooth approximations
$u_\nu\in C^\infty(B_R;\mathbb C^N)
\cap{\rm W}^{1,2}(B_R;\mathbb C^N)$
such that
\[
u_\nu\longrightarrow u
\quad\text{in }{\rm W}^{1,2}(B_R;\mathbb C^N).
\]
Set $F_\nu(\zeta):=T_\zeta u_\nu\in\mathbb X$.
Applying \eqref{e83006} to $u_\nu-u_\mu$, we have
\[
\begin{split}
\|F_\nu-F_\mu\|_{L^2(S^{2n-1};\mathbb X)}^2
&=
\int_{S^{2n-1}}
\|T_\zeta(u_\nu-u_\mu)\|_{\mathbb X}^2
\,d\sigma(\zeta)\\
&\le
C_{n,r,R}
\|u_\nu-u_\mu\|_{{\rm W}^{1,2}(B_R)}^2.
\end{split}
\]
Hence $(F_\nu)$ is Cauchy in
$L^2(S^{2n-1};\mathbb X)$.
Since
$\mathbb X={\rm W}^{1,2}(D_r;\mathbb C^N)$
is a separable Hilbert space, its Bochner $L^2$
space is complete.
Therefore there exists
$F\in L^2(S^{2n-1};\mathbb X)$ such that
\begin{equation}\label{e83008}
F_\nu\longrightarrow F
\quad\text{in }L^2(S^{2n-1};\mathbb X).
\end{equation}

We next extract a subsequence that converges in
${\rm W}^{1,2}(D_r;\mathbb C^N)$
for almost every $\zeta$.
Choose a subsequence $(F_{\nu_k})$ such that
\[
\|F_{\nu_k}-F\|_{L^2(S^{2n-1};\mathbb X)}^2
\le 2^{-k}.
\]
By Tonelli's theorem,
\[
\begin{split}
&\int_{S^{2n-1}}
\sum_{k=1}^{\infty}
\|F_{\nu_k}(\zeta)-F(\zeta)\|_{\mathbb X}^2
\,d\sigma(\zeta)\\
&\qquad=
\sum_{k=1}^{\infty}
\|F_{\nu_k}-F\|_{L^2(S^{2n-1};\mathbb X)}^2
<\infty.
\end{split}
\]
Thus, for almost every $\zeta$,
\[
\sum_{k=1}^{\infty}
\|F_{\nu_k}(\zeta)-F(\zeta)\|_{\mathbb X}^2
<\infty.
\]
In particular,
\[
T_\zeta u_{\nu_k}\longrightarrow F(\zeta)
\quad\text{in }{\rm W}^{1,2}(D_r;\mathbb C^N)
\]
for almost every $\zeta$.

To identify $F(\zeta)$ with $w^{n-1}u(w\zeta)$,
choose a measurable representative of $u$ and define
\[
(Ju)(\zeta,w):=w^{n-1}u(w\zeta).
\]
The complex polar identity \eqref{e82920} gives
\[
\begin{split}
&\|Ju_\nu-Ju\|_{L^2(S^{2n-1}\times D_r)}^2\\
&\qquad=
\int_{S^{2n-1}}\int_{D_r}
|w|^{2n-2}|u_\nu(w\zeta)-u(w\zeta)|^2
\,dA(w)\,d\sigma(\zeta)\\
&\qquad=
2\pi\int_{B_r}|u_\nu(z)-u(z)|^2\,dV(z).
\end{split}
\]
Since $u_\nu\to u$ in $L^2(B_r)$, we obtain
\begin{equation}\label{e83009}
Ju_\nu\longrightarrow Ju
\quad\text{in }L^2(S^{2n-1}\times D_r).
\end{equation}
On the other hand, \eqref{e83008} implies
\begin{equation}\label{e83010}
\begin{split}
F_\nu&\longrightarrow F
\quad\text{in }
L^2\bigl(S^{2n-1};
{\rm W}^{1,2}(D_r;\mathbb C^N)\bigr),\\
&\text{and hence also in }
L^2(S^{2n-1}\times D_r;\mathbb C^N).
\end{split}
\end{equation}
For smooth $u_\nu$,
\[
F_\nu(\zeta,w)
=
Ju_\nu(\zeta,w)
=
w^{n-1}u_\nu(w\zeta).
\]
By uniqueness of the $L^2$ limit,
\eqref{e83009} and \eqref{e83010} yield
\[
F(\zeta,w)=w^{n-1}u(w\zeta)
\]
for almost every $(\zeta,w)$.
Equivalently, for almost every $\zeta$,
\[
F(\zeta)=w^{n-1}u(w\zeta)
\quad\text{in }L^2(D_r;\mathbb C^N).
\]

Since $F(\zeta)\in{\rm W}^{1,2}(D_r;\mathbb C^N)$
for almost every $\zeta$, it follows that
\[
g_\zeta(w):=w^{n-1}u(w\zeta)
\in{\rm W}^{1,2}(D_r;\mathbb C^N)
\]
for almost every $\zeta$.
The estimate in the lemma follows by passing to
the limit in \eqref{e83006}.
This conclusion is independent of the measurable
representative chosen for $u$, because the polar
identity shows that changing $u$ on a null set
changes $Ju$ only on a product-null set.

Finally, passing to the limit in \eqref{e83001}
by the same approximation argument gives,
for almost every $\zeta$,
\begin{equation}\label{e83101}
\bar\partial_w g_\zeta(w)
=
w^{n-1}\sum_{j=1}^n
\bar\zeta_j\bar\partial_{z_j}u(w\zeta)
\quad\text{in }L^2(D_r;\mathbb C^N).
\end{equation}
In particular, \eqref{e83101} holds almost everywhere,
as required.
\end{proof}

\begin{proof}[Proof of Theorem \ref{t01}]
When $n=1$, we have $s=0$, and the conclusion follows
directly from Theorem~\ref{t031} with $v=u$ and $W=V$.
We therefore assume $n\ge2$ throughout the remainder
of the proof.

The idea is to transfer the hypotheses on $u$ to the
compensated slices $g_\zeta$ and then apply the
one-dimensional weighted unique continuation theorem.


It follows from \eqref{e83101} and the condition $|\ddbar u|\leq V|u|$ that
\[
|\bar\partial_w g_\zeta(w)| \le |w|^{n-1} |\bar\partial u(w \zeta)| \le |w|^{n-1} V(w \zeta) |u(w \zeta)| = V(w \zeta) |g_\zeta(w)|.
\]
Set $W_\zeta(w) := |w|^s V(w \zeta)$. Then
\begin{equation}\label{e83102}
|\bar\partial_w g_\zeta| \le |w|^{-s} W_\zeta |g_\zeta|.
\end{equation}
Since
\begin{equation*}
	\begin{split}
\int_{S^{2n-1}} \int_{D_r} |W_\zeta(w)|^{2n} dA(w) d\sigma(\zeta) &= \int_{S^{2n-1}} \int_{D_r} |w|^{2n-2} |V(w \zeta)|^{2n} dA(w) d\sigma(\zeta)\\
&= 2\pi \int_{B_r} |V(z)|^{2n} dV(z) < \infty,
\end{split}
\end{equation*}
$W_\zeta \in L^{2n}(D_r)$ for almost every $\zeta$.

Our goal is to show the sliced $g_\zeta$ is flat at $0$ in the $L^2$ sense as well. Define
\[
G_\zeta(\rho) := \int_{|w|<\rho} |g_\zeta(w)|^2 dA(w).
\]
The polar identity \eqref{e82920} gives
\begin{equation}\label{e83103}
\int_{S^{2n-1}} G_\zeta(\rho) d\sigma(\zeta) = 2\pi \int_{|z|<\rho} |u(z)|^2 dV(z).
\end{equation}

We now use a dyadic argument. Let $\rho_k=2^{-k}$. 
Fix an integer $L\ge1$. Choose $k_L$ sufficiently large
so that $\rho_k<r$ and
\[
\int_{S^{2n-1}}G_\zeta(\rho_k)\,d\sigma(\zeta)
\le \rho_k^{L+2}
\qquad\text{for all }k\ge k_L.
\]
For $k\ge k_L$, set
\[
E_{k,L}:=\{\zeta:G_\zeta(\rho_k)>\rho_k^L\}.
\]
By Markov's inequality,
\[
\sigma(E_{k,L})\le \rho_k^2=4^{-k}.
\]
Consequently,
\[
\sum_{k=k_L}^{\infty}\sigma(E_{k,L})
\le \sum_{k=k_L}^{\infty}4^{-k}<\infty.
\]
The Borel-Cantelli lemma \cite{F} shows immediately that for almost every $\zeta$, there exists $k_0 = k_0(\zeta,L)$ such that for every $k\ge k_0$, $\zeta\notin E_{k,L}$, which means
\begin{equation}\label{e83104}
G_\zeta(\rho_k) \le \rho_k^L \quad \text{for sufficiently large } k.
\end{equation}

For each integer $L\ge 1$, let $Z_L$ be the exceptional null set where \eqref{e83104} fails and $Z=\bigcup_{L=1}^\infty Z_L$. Then $\sigma(Z)=0$. For every $\zeta\notin Z$, and for every integer $L\ge 1$,
\[
G_\zeta(2^{-k}) \le 2^{-kL}
\]
for sufficiently large $k$, where the threshold may depend on $\zeta$ and $L$.

Let $M>0$, and choose an integer $L>M$. Given a small $\rho$, choose $k$ such that $2^{-(k+1)} < \rho \le 2^{-k}$. Since $G_\zeta(\rho)$ is nondecreasing in $\rho$,
\[
G_\zeta(\rho) \le G_\zeta(2^{-k}) \le 2^{-kL}.
\]
Hence
\[
\rho^{-M} G_\zeta(\rho) \le 2^{(k+1)M} 2^{-kL} = 2^M 2^{-k(L-M)} \to 0.
\]
Thus
\begin{equation}\label{e83105}
\lim_{\rho\to 0} \rho^{-M} \int_{|w|<\rho} |g_\zeta(w)|^2 dA(w) = 0
\end{equation}
for any $M>0$ and almost every $\zeta$.

For almost every $\zeta \in S^{2n-1}$, we have
$g_\zeta \in {\rm W}^{1,2}(D_r;\mathbb C^N)$ and
$W_\zeta \in L^{2n}(D_r)$. Moreover, by \eqref{e83102}
and \eqref{e83105}, $g_\zeta$ satisfies the weighted
differential inequality and vanishes to infinite order
at the origin in the $L^2$ sense. Applying
Theorem \ref{t031} on $D_r$, we obtain
\[
g_\zeta = 0 \quad \text{a.e. on } D_r
\]
for almost every $\zeta \in S^{2n-1}$.

The complex polar identity therefore yields
\[
2\pi \int_{B_r} |u(z)|^2\,dV(z)
=
\int_{S^{2n-1}}\int_{D_r}
|g_\zeta(w)|^2\,dA(w)\,d\sigma(\zeta)
=0.
\]
Thus $u=0$ a.e. on $B_r$. Since
$V\in L^{2n}_{\mathrm{loc}}(\Omega)
\subset L^2_{\mathrm{loc}}(\Omega)$ and
$u\in {\rm W}^{1,2}_{\mathrm{loc}}(\Omega;\mathbb C^N)$,
the weak unique continuation theorem of Pan--Zhang
\cite[Theorem~1.3]{PanZhangL2}  applies. As $\Omega$ is connected
and $u$ vanishes on the nonempty open set $B_r$,
we conclude that $u=0$ a.e. on $\Omega$.
This completes the proof.

\end{proof}



We emphasize that the exponent $2n$ is sharp as shown in the following well-known example.
\begin{example}\label{ex1}
Let $0<p<2n$ and $0<\beta<2n/p-1$. Define
\[
u_{\beta}(z):=
\begin{cases}
e^{-|z|^{-\beta}}, & z\neq 0,\\
0, & z=0.
\end{cases}
\]
Then $u_{\beta}$ is smooth, non-trivial and flat at $0$, while
\[
|\bar\partial u_{\beta}| = \frac{\beta}{2}|z|^{-\beta-1}|u_{\beta}|.
\]
Take
\[
V_{\beta}(z)=\frac{\beta}{2}|z|^{-\beta-1},
\]
whose $L^p$-integrability near $0$ is
\[
\int_0^1 r^{2n-1-p(\beta+1)}\,dr < \infty \iff p(\beta+1)<2n.
\]
This works for every $p<2n$, but never for $p=2n$ because then one would require $\beta+1<1$, impossible for $\beta>0$.

\end{example}

 The following example shows that vanishing at a single point
cannot replace the infinite-order vanishing assumption,
even for continuous solutions with
$V\in L^{2n}_{\mathrm{loc}}$.
 
\begin{example}\label{ex2}
Let $B_{e^{-1}}\subset\mathbb C^n$ be the ball centered
at the origin, and define
\[
u(z)=
\begin{cases}
\displaystyle \frac{1}{\log(1/|z|)}, & z\ne0,\\[6pt]
0, & z=0.
\end{cases}
\]
Then $u$ is continuous and nontrivial, with $u(0)=0$.
For $z\ne0$, a direct computation gives
\[
|\nabla u(z)|
=\frac{1}{|z|\log^2(1/|z|)},
\qquad
|\bar\partial u(z)|
=\frac{1}{2|z|\log^2(1/|z|)}.
\]
Consequently,
\[
|\bar\partial u|=V|u|
\quad\text{a.e. on }B_{e^{-1}},
\qquad
V(z):=\frac{1}{2|z|\log(1/|z|)}
\quad(z\ne0),
\]
where we set $V(0)=0$.
The convergence of
\[
\int_0^{e^{-1}}
\frac{dt}{t\log^{4n}(1/t)}
\quad\text{and}\quad
\int_0^{e^{-1}}
\frac{dt}{t\log^{2n}(1/t)}
\]
shows that $u\in W^{1,2n}(B_{e^{-1}})$ and
$V\in L^{2n}(B_{e^{-1}})$, respectively.
However,
\[
\int_{B_r}|u(z)|^2\,dV(z)
\approx \frac{r^{2n}}{\log^2(1/r)}
\qquad\text{as }r\rightarrow 0.
\]
In particular, for every $M>2n$,
\[
r^{-M}\int_{B_r}|u(z)|^2\,dV(z)\longrightarrow+\infty.
\]
Thus $u$ is not $L^2$-flat at the origin. 
\end{example}




\appendix

\refstepcounter{section}
\section*{Appendix A. Failure of a uniform radial Carleman estimate}
\addcontentsline{toc}{section}
  {Appendix A. Failure of a uniform radial Carleman estimate}

In this appendix, $n$ denotes the complex dimension of the source space.
Set
\[
    p_n=\frac{4n}{2n+1},\qquad q_n=\frac{4n}{2n-1},
    \qquad \frac1{p_n}-\frac1{q_n}=\frac1{2n}.
\]
We consider the scale-invariant estimate
\begin{equation}\label{eq:app-model-carleman}
    \bigl\||z|^{-\tau}f\bigr\|_{L^{q_n}(\mathbb C^n)}
    \leq C
    \bigl\||z|^{-\tau}\bar\partial f\bigr\|_{L^{p_n}(\mathbb C^n)},
    \qquad f\in C_c^\infty(\mathbb C^n\setminus\{0\}),
\end{equation}
with a constant independent of the large parameter $\tau$.
Here $dV$ is standard Lebesgue measure and
\[
    |\bar\partial f|^2
    :=\sum_{j=1}^n\left|\frac{\partial f}{\partial\bar z_j}\right|^2.
\]
We show that such a uniform estimate fails when $n\geq2$, even
for test functions supported in a fixed compact annulus.

\begin{proposition}\label{prop:app-carleman-failure}
Let $n\geq2$. There exist constants $c,C,\tau_0>0$ and nonzero
functions $f_\tau\in C_c^\infty(\mathbb C^n\setminus\{0\})$,
defined for every real $\tau\geq\tau_0$, such that
\[
    \operatorname{supp}f_\tau
    \subset\{e^{-1/2}\leq|z|\leq e^{3/2}\}
\]
and
\begin{equation}\label{eq:app-ratio-growth}
    c\tau^{(n-1)/(2n)}
    \leq
    \frac{\bigl\||z|^{-\tau}f_\tau\bigr\|_{L^{q_n}}}
         {\bigl\||z|^{-\tau}\bar\partial f_\tau\bigr\|_{L^{p_n}}}
    \leq C\tau^{(n-1)/(2n)}.
\end{equation}
Consequently, \eqref{eq:app-model-carleman} cannot hold with a
uniform constant along any unbounded set of positive parameters.
\end{proposition}

The construction uses holomorphic monomials whose restrictions
to the unit sphere concentrate near $\{|\omega_1|=1\}$.
We first record their exact spherical moments.

\begin{lemma}\label{lem:app-sphere-moments}
Let $d\sigma$ denote unnormalized Euclidean surface measure on
$S^{2n-1}$. For $n\geq1$ and $s>-1$,
\begin{equation}\label{eq:app-sphere-moment}
    M_n(s):=\int_{S^{2n-1}}|\omega_1|^{2s}\,d\sigma(\omega)
    =\frac{2\pi^n\Gamma(s+1)}{\Gamma(s+n)}.
\end{equation}
In particular, for every fixed $0<\ell<\infty$,
\begin{equation}\label{eq:app-angular-norm}
    \|\omega_1^m\|_{L^\ell(S^{2n-1})}
    \approx m^{-(n-1)/\ell},\qquad m\to\infty,
\end{equation}
where the comparison constants depend only on $n$ and $\ell$.
\end{lemma}

\begin{proof}
Computing a Gaussian integral in Cartesian coordinates gives
\[
    \int_{\mathbb C^n}|z_1|^{2s}e^{-|z|^2}\,dV(z)
    =\left(2\pi\int_0^\infty\rho^{2s+1}e^{-\rho^2}\,d\rho\right)
      \pi^{n-1}
    =\pi^n\Gamma(s+1).
\]
On the other hand, polar coordinates yield
\[
    \int_{\mathbb C^n}|z_1|^{2s}e^{-|z|^2}\,dV(z)
    =M_n(s)\int_0^\infty r^{2s+2n-1}e^{-r^2}\,dr
    =\frac12M_n(s)\Gamma(s+n).
\]
Comparing the two expressions proves \eqref{eq:app-sphere-moment}.
For $n\geq2$,
\[
    \frac{\Gamma(s+1)}{\Gamma(s+n)}
    =\frac1{(s+1)\cdots(s+n-1)}\approx s^{-(n-1)}
    \qquad(s\to\infty).
\]
Since
$\|\omega_1^m\|_{L^\ell}^{\ell}=M_n(m\ell/2)$,
this proves \eqref{eq:app-angular-norm}.
For $n=1$, the conclusion follows directly from $|\omega_1|=1$
on $S^1$.
\end{proof}

For comparison with the usual beta-integral formula, the same
surface-measure convention gives, for $n\geq2$,
\[
    \frac1{\sigma(S^{2n-1})}
    \int_{S^{2n-1}}F(|\omega_1|^2)\,d\sigma(\omega)
    =(n-1)\int_0^1F(u)(1-u)^{n-2}\,du,
\]
where $\sigma(S^{2n-1})=2\pi^n/\Gamma(n)$.

\begin{proof}[Proof of Proposition~\ref{prop:app-carleman-failure}]
Fix a real-valued function $\chi\in C_c^\infty(\mathbb R)$ such that
\[
    0\leq\chi\leq1,\qquad
    \operatorname{supp}\chi\subset[-\tfrac12,\tfrac32],\qquad
    \chi=1\quad\text{on }[0,1].
\]
For every sufficiently large $\tau$, choose
\[
    m=\left\lfloor\tau-n+\frac12\right\rfloor,
    \qquad \beta=\tau-m-n+\frac12.
\]
Then $m\geq1$, $0\leq\beta<1$, and $m=\tau+O_n(1)$.
Define
\[
    f_\tau(z)=z_1^m|z|^\beta\chi(\log|z|).
\]
The support lies in the fixed annulus appearing in the proposition,
so $f_\tau$ is smooth and compactly supported away from the origin,
including when $\beta$ is nonintegral.

Write $p=p_n$, $q=q_n$, and $z=r\omega$, where $r=|z|$ and
$\omega\in S^{2n-1}$. The choice of $\beta$ gives the two identities
\begin{equation}\label{eq:app-radial-balance}
    m+\beta-\tau=-n+\frac12=-\frac{2n}{q},\qquad
    m+\beta-\tau-1=-n-\frac12=-\frac{2n}{p}.
\end{equation}
Using $dV=r^{2n-1}\,dr\,d\sigma$ and then $t=\log r$, we obtain
\begin{equation}\label{eq:app-left-norm}
\begin{aligned}
    \bigl\||z|^{-\tau}f_\tau\bigr\|_{L^q}^q
    &=M_n(mq/2)\int_0^\infty
      r^{(m+\beta-\tau)q+2n-1}|\chi(\log r)|^q\,dr\\
    &=M_n(mq/2)\int_0^\infty|\chi(\log r)|^q\,\frac{dr}{r}\\
    &=M_n(mq/2)\|\chi\|_{L^q(\mathbb R)}^q.
\end{aligned}
\end{equation}

To compute the derivative, observe that $\bar\partial(z_1^m)=0$
and
\[
    \frac{\partial r}{\partial\bar z_j}=\frac{z_j}{2r},
    \qquad |\bar\partial r|=\frac12.
\]
It follows that
\begin{equation}\label{eq:app-dbar-pointwise}
    |\bar\partial f_\tau(r\omega)|
    =\frac12r^{m+\beta-1}|\omega_1|^m
       |\chi'(\log r)+\beta\chi(\log r)|.
\end{equation}
In particular, differentiating the holomorphic monomial produces
no factor of $m$. Applying the second identity in
\eqref{eq:app-radial-balance}, we find
\begin{equation}\label{eq:app-right-norm}
\begin{aligned}
    \bigl\||z|^{-\tau}\bar\partial f_\tau\bigr\|_{L^p}^p
    &=2^{-p}M_n(mp/2)\int_0^\infty
      r^{(m+\beta-\tau-1)p+2n-1}
      |\chi'(\log r)+\beta\chi(\log r)|^p\,dr\\
    &=2^{-p}M_n(mp/2)
      \|\chi'+\beta\chi\|_{L^p(\mathbb R)}^p.
\end{aligned}
\end{equation}

The radial factor is bounded above and below uniformly for
$\beta\in[0,1]$. Indeed, the map
$\beta\mapsto\|\chi'+\beta\chi\|_{L^p}$ is continuous.
If it vanished for some $\beta$, then
$(e^{\beta t}\chi(t))'=0$. Compact support would force
$\chi\equiv0$, a contradiction. Compactness of $[0,1]$ therefore
gives constants $c_0,C_0>0$ such that
\[
    c_0\leq\|\chi'+\beta\chi\|_{L^p}\leq C_0,
    \qquad 0\leq\beta\leq1.
\]

Combining \eqref{eq:app-left-norm}, \eqref{eq:app-right-norm},
and Lemma~\ref{lem:app-sphere-moments}, we obtain
\[
\begin{aligned}
    \frac{\bigl\||z|^{-\tau}f_\tau\bigr\|_{L^q}}
         {\bigl\||z|^{-\tau}\bar\partial f_\tau\bigr\|_{L^p}}
    &=2\frac{\|\chi\|_{L^q}}{\|\chi'+\beta\chi\|_{L^p}}
       \frac{M_n(mq/2)^{1/q}}{M_n(mp/2)^{1/p}}\\
    &\approx m^{(n-1)(1/p-1/q)}
     =m^{(n-1)/(2n)}
     \approx\tau^{(n-1)/(2n)}.
\end{aligned}
\]
All comparison constants are independent of $\tau$, $m$, and
$\beta$. Since $n\geq2$, the last quantity tends to infinity,
which proves the proposition.
\end{proof}

The proposition gives a necessary lower bound
$C_\tau\geq c\tau^{(n-1)/(2n)}$ for any constants in
\eqref{eq:app-model-carleman}. It establishes failure of uniformity
in $\tau$; it does not by itself rule out estimates whose constants
depend on $\tau$. The obstruction comes from angular concentration
of holomorphic monomials, which is not penalized by $\bar\partial$.


When $n=1$, a uniform estimate does hold for integer parameters.
The Cauchy--Pompeiu formula and the Hardy--Littlewood--Sobolev
inequality, equivalently the unweighted case of the Stein--Weiss
inequality \cite{SW}, give
\[
    \|g\|_{L^4(\mathbb C)}
    \leq C\|\bar\partial g\|_{L^{4/3}(\mathbb C)},
    \qquad g\in C_c^\infty(\mathbb C).
\]
For $f\in C_c^\infty(\mathbb C\setminus\{0\})$ and $k\in\mathbb Z$,
apply this estimate to $g=z^{-k}f$, extended by zero near the origin.
Since $\bar\partial g=z^{-k}\bar\partial f$, this gives
\begin{equation}\label{eq:app-planar-integers}
    \bigl\||z|^{-k}f\bigr\|_{L^4(\mathbb C)}
    \leq C\bigl\||z|^{-k}\bar\partial f\bigr\|_{L^{4/3}(\mathbb C)},
\end{equation}
with a constant independent of $k$.

The parameter restriction matters: at each half-integer
$\tau=k+\tfrac12$, the estimate on all of
$C_c^\infty(\mathbb C\setminus\{0\})$ fails even for this fixed
$\tau$. To see this, fix $k\in\mathbb Z$ and, for $L\geq1$, set
\[
    F_L(z)=z^k\chi\!\left(\frac{\log|z|}{L}\right).
\]
Each $F_L$ is smooth and compactly supported away from the origin.
Polar coordinates give
\[
    \bigl\||z|^{-k-1/2}F_L\bigr\|_{L^4}
    =(2\pi)^{1/4}L^{1/4}\|\chi\|_{L^4},
\]
whereas
\[
    \bigl\||z|^{-k-1/2}\bar\partial F_L\bigr\|_{L^{4/3}}
    =\frac12(2\pi)^{3/4}L^{-1/4}\|\chi'\|_{L^{4/3}}.
\]
Their ratio is a positive constant times $L^{1/2}$ and hence
diverges as $L\to\infty$. Thus the planar assertion should be
stated with an admissible parameter set, such as the integers
in \eqref{eq:app-planar-integers}.



\refstepcounter{section}
\section*{Appendix B. A local Hardy inequality}
\addcontentsline{toc}{section}
  {Appendix B. A local Hardy inequality}
Throughout this appendix, we assume $n\ge2$ and set
$d=2n\ge4$. Identifying $\mathbb C^n$ with $\mathbb R^d$,
we prove the local Hardy inequality for
$h\in {\rm W}^{1,2} (B_R;\mathbb C^N)$ used in the compensated
slicing argument.

Recall that the classical Hardy inequality gives
\begin{equation}\label{e90201}
\frac{(d-2)^2}{4}
\int_{\mathbb R^d}\frac{|f(x)|^2}{|x|^2}\,dx
\le
\int_{\mathbb R^d}|\nabla f(x)|^2\,dx,
\qquad
f\in {\rm W}^{1,2}(\mathbb R^d;\mathbb C^N).
\end{equation}


Recall that $0<r<R$. More explicitly, we will show the following
\begin{equation}\label{e90202}
\int_{B_r}\frac{|h(z)|^2}{|z|^2}dV(z) \le \frac{2}{(n-1)^2}\left(\int_{B_R}|\nabla h(z)|^2dV(z) + \frac{16}{(R-r)^2}\int_{B_R}|h(z)|^2dV(z)\right).
\end{equation}

\begin{proof}
Choose a smooth function $\theta:\mathbb{R}\to[0,1]$ such that
\[
\theta(t)=1 \quad \text{for} \quad t\le 0,\quad \theta(t)=0 \quad \text{for} \quad t\ge 1,\quad |\theta'(t)|\le 2.
\]
Define
\[
\eta(x):=\theta\left(\frac{|x|-r}{R_1-r}\right),
\]
where $R_1=(r+R)/2$. Then $\eta$ is smooth, $0\le \eta\le 1$, $\eta=1$ on $B_r$, and $\operatorname{supp}\eta\subset \bar B_{R_1}\subset B_R$. Moreover, we have in $r<|x|<R_1$,
\[
|\nabla \eta(x)|\le \frac{2}{R_1-r} = \frac{4}{R-r}.
\]

Recall that $h\in {\rm W}^{1,2}(B_R;\mathbb{C}^N)$. Set $f=\eta h$. Then $f\in {\rm W}^{1,2}(B_R)$ with $\operatorname{supp} f\subset \overline{B}_{R_1}\Subset B_R$. We extend $f$ by zero outside $B_R$ and denote the extension by $\tilde f\in {\rm W}^{1,2}(\mathbb{R}^{2n})$.
Then
\begin{equation}\label{e90203}
\int_{B_r}\frac{|h(z)|^2}{|z|^2}dV(z) = \int_{B_r}\frac{|\eta h(x)|^2}{|x|^2}\,dx \le \int_{\mathbb{R}^d}\frac{|\tilde f(x)|^2}{|x|^2}\,dx.
\end{equation}
We deduce from the global Hardy inequality \eqref{e90201} that
\begin{equation}\label{e90204}
\begin{split}
\int_{\mathbb{R}^d}\frac{|\tilde f(x)|^2}{|x|^2}\,dx &\le \frac{1}{(n-1)^2}\int_{\mathbb{R}^d}|\nabla \tilde f(x)|^2\,dx= \frac{1}{(n-1)^2}\int_{B_R}|\eta \nabla h + h\nabla \eta|^2dx\\
&\le \frac{2}{(n-1)^2}\int_{B_R}\eta^2|\nabla h|^2dx + \frac{2}{(n-1)^2}\int_{B_R}|h|^2|\nabla \eta|^2dx \\
&\le \frac{2}{(n-1)^2}\int_{B_R}|\nabla h|^2dx + \frac{2}{(n-1)^2}\|\nabla \eta\|_{L^\infty}^2\int_{B_R}|h|^2dx \\
&\le \frac{2}{(n-1)^2}\left(\int_{B_R}|\nabla h|^2dV(z) + \frac{16}{(R-r)^2}\int_{B_R}|h|^2dV(z)\right).
\end{split}
\end{equation}
Combining \eqref{e90203} and \eqref{e90204} gives the local Hardy inequality \eqref{e90202}.
\end{proof}

\clearpage
\phantomsection
\addcontentsline{toc}{section}{References}


\begin{thebibliography}{99}

\bibitem{A} N. Aronszajn, \emph{A unique continuation theorem for solutions of elliptic partial differential equations or inequalities of second order},  J. Math. Pures Appl., (9) \textbf{36} (1957), 235-249.



\bibitem{BellLempert} S. Bell and L. Lempert, \emph{A $C^\infty$ Schwarz reflection principle in one and several complex variables}, J. Differential Geom.,\textbf{32} (1990), no.3, 899-915.

\bibitem{C} T. Carleman, \emph{Sur un problème d'unicité pour les systèmes d'équations aux dérivées partielles à deux variables indépendantes}, Ark. Mat. Astr. Fys., 26B (1939), no. 17, 1–9. 


\bibitem{F} W. Feller, \emph{An introduction to probability theory and its applications}, Vol. 1 (3rd ed.). New York: Wiley, (1968).


\bibitem{H} L. Hörmander, \emph{The analysis of linear partial differential operators III: pseudo-differential operators}, Springer Berlin, Heidelberg (2007).

\bibitem{JK} D. Jerison and C. E. Kenig, \emph{Unique continuation and absence of positive eigenvalues for Schr\"odinger operators. With an appendix by E. M. Stein}, Ann. of Math., (2)  \textbf{121} (1985), no.3, 463-494.


\bibitem{Mandache} N. Mandache, \emph{A counterexample to unique continuation in dimension two}, Comm. Anal. Geom., \textbf{10} (2002), no.1, 1-10.



\bibitem{PZ} Y. Pan and Y. Zhang, \emph{Unique continuation for $\ddbar u=Vu$ at infinity}, arXiv:2608.22698.


\bibitem{PZ2025} Y. Pan and Y. Zhang, \emph{Unique continuation of Schrödinger-type equations for $\bar\partial$}, Comm. Anal. Geom., \textbf{33} (2025), no. 8, 1975-2014.


\bibitem{PanZhang2} Y. Pan and Y. Zhang, \emph{A remark on unique continuation of $\bar{\partial}$}, Nagoya Math. J., \textbf{261}  (2026), e55 1–12.



\bibitem{PanZhangL2} Y. Pan and Y. Zhang, \emph{Unique continuation for $\bar{\partial}$ with square-integrable potentials}, New York J. Math., \textbf{29} (2023), 402-416.


\bibitem{Shi} Z. Shi, \emph{A unique continuation property for $|\bar{\partial} u| \leq V|u|$}, Math. Scand., \textbf{131} (2025), no.1, 113-154.


\bibitem{SW} E.-M. Stein and G. Weiss, \emph{Fractional integrals on $n$-dimensional Euclidean space}, J. Math. Mech., \textbf{7} (1958), 503-514.
	
\end{thebibliography}
\end{document}